\documentclass[preprint,3p,a4paper]{elsarticle}
\usepackage{mathrsfs}
\usepackage{amsfonts,amsmath,amssymb,amsthm}

\usepackage{latexsym}

\usepackage{hyperref}
\usepackage{geometry}
 \theoremstyle{plain}
  \newtheorem{thm}{Theorem}[section]
  \newtheorem{lem}[thm]{Lemma}
  \newtheorem{prop}[thm]{Proposition}
  \newtheorem{cor}[thm]{Corollary}
  \newtheorem{defn}[thm]{Definition}
  
  \newtheorem{exam}[thm]{Example}
  \newtheorem{rem}[thm]{Remark}
  
  \newtheorem{fact}[thm]{Fact}

\newcommand{\da}{\downarrow\hspace{-2pt}}

\newcommand{\ua}{\uparrow\hspace{-2pt}}

\newcommand{\ra}{\rightarrow}

\newcommand{\UUa}{\mathord{\mbox{\makebox[0pt][l]{\raisebox{.4ex}{$\uparrow$}}$\uparrow$}}}

\journal{Topology and its applications}
\UseRawInputEncoding
\begin{document}

\begin{frontmatter}




\title{Two topologise  on  the lattice of  Scott closed subsets\tnoteref{t1}}
\tnotetext[t1]{Research supported by NSF of China (Nos. 11871353, 12001385).}

\author{Yu Chen}
\ead{eidolon\_chenyu@126.com}
\author{Hui Kou\corref{cor}}
\ead{kouhui@scu.edu.cn}
\author{Zhenchao Lyu\corref{}}
\ead{zhenchaolyu@scu.edu.cn}
\cortext[cor]{Corresponding author}
\address{Department of Mathematics, Sichuan University, Chengdu 610064, China}

\begin{abstract}
\indent For a poset $P$, let $\sigma(P)$ and $\Gamma(P)$ respectively denote the lattice of its Scott open subsets and Scott closed subsets ordered by inclusion, and set $\Sigma P=(P,\sigma(P))$.  In this paper, we discuss the lower Vietoris topology and the Scott topology on $\Gamma(P)$ and give some sufficient conditions to make the two topologies equal. We built an adjunction between $\sigma(P)$ and $\sigma(\Gamma(P))$ and proved that
$\Sigma P$ is core-compact iff $\Sigma\Gamma(P)$ is core-compact iff $\Sigma\Gamma(P)$ is sober, locally compact and $\sigma(\Gamma(P))=\upsilon(\Gamma(P))$ (the lower Vietoris topology). This answers a question in \cite{Miao}.
Brecht and Kawai \cite{Matthew} asked whether the consonance of a topological space $X$ implies the consonance of its lower powerspace, we give a partial answer to this question at the last part of this paper.

\end{abstract}

\end{frontmatter}

\section{Introduction}
The lower powerspace over a topological space $X$ is the set of closed subsets of $X$ with the lower Vietoris topology. The lower powerspace coincides with the lower powerdomain for continuous dcpos with the Scott topology, where the latter construction is used in modelling non-deterministic computation (see \cite{A,smyth}). There naturally arise a question when the lower Vietoris topology and the Scott topology coincide on the lattice of closed subsets of a topological space. For a poset $P$, we give some sufficient conditions such that the lower Vietoris topology and the Scott topology coincide on $\Gamma(P)$. We observe that there is an adjunction between $\sigma(P)$ and $\sigma(\Gamma(P))$, which serves as a useful tool in studying the relation between $P$ and $\Gamma(P)$. Particularly, we obtain the conclusion that $(P, \sigma(P))$ is core-compact iff $(\Gamma(P), \sigma(\Gamma(P))$ is a sober and locally compact space. This positively answers a question in \cite{Miao}. We show that the two topologies generally don't coincide on the lattice of closed subsets of a topological space that is not a Scott space.

In domain theory, the spectral theory of distributive continuous lattices and the duality theorem have been extensively studied. Hoffmann \cite{ReHof} and Lawson \cite{Law} independently proved that a directed complete poset (or dcpo for short) is continuous iff the lattice of its Scott open subsets is completely distributive. Gierz and Lawson \cite {Gierz} proved that quasicontinuous dcpos equipped with the Scott topologies are precisely the spectrum of distributive hypercontinuous lattices. The correspondence between $P$ and $\Gamma(P)$ has also been studied (see \cite{redbook,Hof,Law,ven,HoZhao}). A poset $P$ is continuous (resp., quasicontinuous, algebraic, quasialgebraic) if and only if $\Gamma(P) $ is continuous (resp., quasicontinuous, algebraic, quasialgebraic). Using the pair of adjoint we can give a unified proof.

We see that if $L$ is a continuous lattice such that $\sigma(L^{op})=\upsilon(L^{op})$, then $(L^{op}, \sigma(L^{op}))$ is a sober and locally compact space. We give two equivalent conditions of $\sigma(L^{op})=\upsilon(L^{op})$ for continuous lattices.
The previous result about core-compact posets implies that if $L$ is a distributive continuous lattice such that the hull-kernel topology of Spec$L$ is just the Scott topology then $\sigma(L^{op})=\upsilon(L^{op})$.
We wonder whether the condition $\sigma(L^{op})=\upsilon(L^{op})$ implies that Spec$L$ is a Scott space. If not, what other conditions are needed. Actually, characterize those distributive continuous lattices for which the spectrum is a dcpo equipped with the Scott
topology is an open problem \cite{Van}.
The case of distributive algebraic lattices is simpler. The notion of jointly Scott continuity occurs in our proofs. It is well-known that a complete lattice is sober with respect to the Scott topology if it is  jointly Scott continuity \cite[Corollary II-1.12]{redbook}.  We give an example of a complete lattice   which is sober with  the Scott topology but fails to be  jointly Scott continuity.

Consonance is an important topological property (see \cite{Matthew,Bouziad,Dolecki,Nogura}). It is proved in \cite{Dolecki} that every regular C\'{e}ch-complete space is consonant. Thus every completely metrizable space is consonant. Bouziad \cite{Bouziad} proved that the set of rational numbers with subspace topology from $\mathbb R$ is not consonant. A space $X$ is consonant if and only if, for every topological space $Y$, the compact-open and the Isbell topologies agree on the space of continuous maps from $X$ to $Y$ (see \cite{Jean-2013}). Brecht and Kawai \cite{Matthew} proved that the consonance of a $T_0$ space $X$ is equivalent to the commutativity of the upper and lower powerspaces of $X$. There they asked whether the consonance of $X$ implies the consonance of its lower powerspace. We will give a partial answer at the last part of this paper.

\section{Preliminaries}
We refer to (\cite{A,redbook,Jean-2013}) for some concepts and notations of domain theory that will be used in the paper.

For a poset $P$ and $A\subseteq P$, let
$\da A=\{x\in P: \exists a\in A, \ x\leq a\}$ and $\ua A=\{x\in P: \exists a\in P, \ a\leq x\}$. For  $x\in P$, we write
$\da x$ for $\da\{x\}$ and
$\ua x$ for $\ua\{x\}$. A set $A$ is called saturated if and only if $A=\ua A$.
If $F$ is a finite subset of $P$, we denote by $F\subseteq_f\!P$.
A poset $P$ is called a directed complete poset (dcpo, for short) if every directed subset $D$ of $P$ has a supremum.

For a poset $P$, the upper topology $\upsilon(P)$ is the topology generated by taking the collection of sets $\{P\setminus\da x: x\in P\}$ as a subbais, the lower topology $\omega(P)$ on $P$ is defined dually. A subset $A$ of $P$ is called Scott closed
if $A=\da A$ and for any directed set $D\subseteq A$, $\bigvee D\in A$ whenever the least upper bound $\bigvee D$ exists. The Scott topology $\sigma(P)$ consists of the complements of all Scott closed sets of $P$. The topology $\lambda(P)=\sigma(P)\vee\omega(P)$ is called the Lawson topology on $P$.

Given a topological space $X$, let $\mathcal O(X)$ be the topology on $X$, then $\mathcal O(X)$ is a complete lattice  ordered by inclusion. We denote by $(C(X),\upsilon(C(X)))$ the lattice of closed subsets of $X$ equipped with the lower Vietoris topology, which is generated by $\Diamond U =\{A\in C(X): A\cap U\ne\emptyset\}$ as a subbase, where $U$ ranges over $\mathcal O(X)$. A topological space is core-compact if and only if $\mathcal O(X)$ is a continuous lattice. An arbitrary nonempty subset $A$ of a topological space $X$ is irreducible if $A\subseteq B\cup C$ for closed subsets $B$ and $C$ implies $A\subseteq B$ or $A\subseteq C$. $X$ is called sober if for every irreducible closed set $A$, there exists a unique $x\in X$ such that $\{x\}^-=A$.

For a poset $P$ and $x,y\in P$, we say that $x$ is way-below $y$, in symbols $x\ll y$, if for any directed subset $D$ of $P$ that has a least upper bound in $P$, $y\leq\bigvee D$ implies $x\leq d$ for some $d\in D$. For any non-empty subsets $F, G$ of $P$, we say that $F$ is way-below $G$ and write $F\ll G$ if for any directed subset $D$ of $P$ that has a least upper bound in $P$, $\bigvee D\in\ua G$ implies $D\,\cap\ua F\neq\emptyset$. $P$ is continuous if the set $\{d\in D: d\ll x\}$ is directed and $x=\bigvee\{d\in P: d\ll x\}$ for all $x\in P$. Let $\text{K}(P)=\{x\in P: x\ll x\}$. $P$ is algebraic if for all $x\in P$ the set $\da x\cap\text K(P)$ is directed and $x=\bigvee(\da x\cap\text K(P))$.

A poset $P$ is called a quasicontinuous poset (resp., quasialgebraic poset) if for all $x\in P$ and $U\in\sigma(P), x\in U$ implies that there is a non-empty finite set $F\subseteq P$ such that $x\in\text{int}_{\sigma}\ua F\subseteq\ua F\subseteq U$ (resp., $x\in\text{int}_{\sigma}\ua F=\ua F\subseteq U$).

For a complete lattice $L$ and $x,y\in L$, we say that $x\prec y$ iff whenever the intersection of a non-empty collection of upper sets is contained in $\ua y$, then the intersection of finitely many is contained in $\ua x$. $x\prec y\Leftrightarrow y\in\text{int}_{\upsilon}\ua x$ and $x\prec y$ implies $x\ll y$. $L$ is hypercontinuous if $\{d\in L: d\prec x\}$ is directed and $x=\bigvee\{d\in L: d\prec x\}$ for all $x\in L$. $L$ is hyperalgebraic provided $x=\bigvee\{y\in L: y\prec y\leq x\}$ for all $x\in L$.

Let $L$ be a complete lattice, we say that $x$ is wayway-below $y$, in notation $x\lll y$, if for any subset $A$ of $L$, $y\leq\bigvee A$ implies $x\leq a$ for some $a\in A$. $L$ is  prime continuous  if $x=\bigvee\{d\in L: d\lll x\}$ for all $x\in L$.

The next two well-known propositions give various equivalent formulations for completely distributive lattice and hypercontinuous lattice separately.

\begin{prop}\rm\label{CDlattice}(\cite{redbook})
Let $L$ be a complete lattice, the following conditions are equivalent:
\begin{enumerate}
\item[(1)] $L$ is prime-continuous,
\item[(2)] $L$ is a completely distributive lattice,
\item[(3)] $L$ is distributive and both $L$ and $L^{op}$ are continuous lattices,
\end{enumerate}
\end{prop}

\begin{thm}\rm\label{property of hypercontinuous}(\cite{Gierz})
Let $L$ be a complete lattice, then the following conditions are equivalent:
\begin{enumerate}
\item[(1)] $L$ is a hypercontinuous lattice,
\item[(2)] $L$ is a continuous lattice and $\sigma(L)=\upsilon(L)$,
\item[(3)] $L$ is a continuous lattice, $L^{op}$ is a quasicontinuous lattice, and the bi-Scott topology agrees with the Lawson topology$(\lambda(L)=\sigma(L)\vee\sigma(L^{op}))$.
\end{enumerate}
\end{thm}

The following lemma characterizes the binary relation $\prec$ on the complete lattice of Scott open subsets. It is easy to see that this conclusion is also true for any topological space.
\begin{lem}\rm\label{approximate of hypercontinuous}
Let $P$ be a poset and $U, V\in\sigma(P)$. $U\prec V$ if and only if $U\subseteq\ua F\subseteq V$ for some $F\subseteq_f\!P$.
\end{lem}

\noindent{\bf Proof.} $(\Leftarrow)$ Let $F=\{x_1, x_2,\cdots, x_n\}$, then $V\nsubseteq P\setminus\da x_i$, for each $1\leq i\leq n$. We claim that $V\in\bigcap_{i=1}^n(\sigma(P)\setminus\da_{\sigma(P)}(P\setminus\da x_i))\subseteq\ua_{\sigma(P)}U$. Indeed, for each $W\in\bigcap_{i=1}^n(\sigma(P)\setminus\da_{\sigma(P)}(P\setminus\da x_i))$, $W\cap\da x_i\neq\emptyset$ for each $1\leq i\leq n$, then $U\subseteq\ua F\subseteq W$.

$(\Rightarrow)$ If $U\prec V$, then there exists $\{W_i\in\sigma(P): 1\leq i\leq n\}$ such that $V\in\bigcap_{i=1}^n(\sigma(P)\setminus\da_{\sigma(P)}W_i)\subseteq\ua_{\sigma(P)}U$, which implies that $V\nsubseteq W_i$ for each $i$. Let $F=\{x_1, x_2,\cdots, x_n\}$, where $x_i\in V\setminus W_i$ for each $1\leq i\leq n$. Now we show that $U\subseteq\ua F$. If not, there exists $x\in U\setminus\ua F$, i.e., $F\subseteq P\setminus\da x$. Then $P\setminus\da x\nsubseteq W_i$ for each $1\leq i\leq n$, but $U\nsubseteq P\setminus\da x$, which is a contradiction. $\hfill{} \square$
\vskip 3mm
M. Ern\'e \cite{e2009} and T. Yokoyama \cite{Yokoyama} respectively give the following Lemma, and prove the result that a spectral Scott space (i.e. a poset on which the Scott topology is spectral) is a quasialgebraic domain.

\begin{lem}\rm\label{Scott compact open}
Let $P$ be a poset. A Scott open subset $U$ of $P$ is Scott compact if and only if there exists $F\subseteq_f\!P$ such that $U=\text{int}_{\sigma}\ua F=\ua F$.
\end{lem}

The correspondence between $P$ and $\sigma(P)$ is enumerated below, where (2) and (3) is a direct consequences of Lemma \ref{approximate of hypercontinuous} and \ref{Scott compact open}.

\begin{thm}\rm\label{dual property of opensubsets}(\cite{redbook,mao,xul})
Let $P$ be a poset,
\begin{enumerate}
\item[(1)] $P$ is continuous (algebraic) iff the lattice $\sigma(P)$ of all Scott open sets is a completely distributive lattice (completely distributive algebraic lattice).
\item[(2)] $P$ is quasicontinuous iff the lattice $\sigma(P)$ of all Scott open sets is a hypercontinuous lattice.
\item[(3)] $P$ is quasialgebraic iff the lattice $\sigma(P)$ of all Scott open sets is an algebraic lattice iff $\sigma(P)$ is a hyperalgebraic lattice.
\end{enumerate}
\end{thm}

\section{The correspondent properties between $P$ and $\Gamma(P)$}
Recently, H. Miao et al \cite{Miao} proved that for a well-filtered dcpo $L$, $\Sigma L$ is locally compact if and only if  $\Sigma\Gamma(L)$ is a  locally compact space. They further asked the following question \cite[Problem 5.9]{Miao}:

For a poset $P$, if $\Sigma P$ is core-compact, must $\Sigma\Gamma(P)$ be locally compact?

In this section, we investigate conditions when the lower Vietoris topology and the Scott topology on $\Gamma(P)$ coincide for a poset $P$ and give an positive answer to the above question. Particularly, we show that for a poset $L$, $\Sigma P$ is core-compact iff  $\Sigma\Gamma(P)$ is core-compact iff $\Sigma\Gamma(P)$ is sober, locally compact and $\sigma(\Gamma(P))=\upsilon(\Gamma(P))$ .

When we are working on continuous lattices, we notice that the condition $\sigma(L^{op})=\upsilon(L^{op})$ is very important. Firstly, several well-known lemmas are needed.

\begin{lem}\rm\label{priestly space}(\cite{Gierz,ven})
A complete lattice $L$ is a quasicontinuous (quasialgebraic) lattice iff $\omega(L)$ is a continuous (algebraic) lattice.
\end{lem}

\begin{lem}\rm\label{sober corecompact lc}(\cite[Theorem V-5.6.]{redbook})
For a sober space $X$, the lattice $\mathcal O(X)$ of open subsets is
continuous iff $X$ is locally compact.
\end{lem}

The following lemma is easy to verify.

\begin{lem}\rm\label{v is sober}
For a complete lattice $L$, both $(L, \upsilon(L))$ and $(L, \omega(L))$ are sober.
\end{lem}

\begin{prop}\rm\label{Lop is sober locally compact}
Let $L$ be a continuous lattice, if $L$ satisfies the condition that $\upsilon(L^{op})=\sigma(L^{op})$, then $(L^{op}, \sigma(L^{op}))$ is a sober and locally compact space.
\end{prop}
\noindent{\bf Proof.}
By Lemma \ref{priestly space}, $\sigma(L^{op})=\upsilon(L^{op})=\omega(L)$ are continuous lattices. Thus $(L, \sigma(L^{op}))$ is a sober and locally compact space by Lemma \ref{sober corecompact lc} and \ref{v is sober} .
$\hfill{} \square$
\vskip 3mm
We give a sufficient condition such that $\sigma(\Gamma(P))=\upsilon(\Gamma(P))$ for a poset $P$, which is crucial for further discussion.

\begin{prop}\rm\label{Scott=upper}
Let $P$ be a poset. If $\Sigma(\prod\limits^nP)=\prod\limits^n(\Sigma P)$ for each $n\in\mathbb N$, then $\sigma(\Gamma(P))=\upsilon(\Gamma(P))$.
\end{prop}

\noindent{\bf Proof.}
We only need to prove that $\sigma(\Gamma(P))\subseteq\upsilon(\Gamma(P))$.

At first, for each $n\in\mathbb{N}$, we define a map $s_n: \prod\limits^n P\rightarrow\Gamma(P)$ as follows:
 $$\textstyle\forall  (x_1,x_2,\ldots,x_n) \in \prod\limits^nP, \displaystyle\ s_n(x_1,x_2,\ldots,x_n)=\bigcup\limits_{k=1}^n\downarrow\!x_k.$$
  We claim that $s_n$ preseves existing directed sups. Let $\{(x_{1i},x_{2i},\ldots,x_{ni}): \ i\in I\}$ be a directed subset of  $\prod\limits^nP$ such that the supremum of $\{(x_{1i},x_{2i},\ldots,x_{ni}): \ i\in I\}$ exists in $\prod\limits^nP$, then  $\bigvee_{i\in I}(x_{1i},x_{2i},\ldots,x_{ni})=(\bigvee_{i\in I}x_{1i},\bigvee_{i\in I}x_{2i},\ldots,\bigvee_{i\in I}x_{ni})$. We have
$$
s_n(\bigvee_{i\in I}(x_{1i},x_{2i},\ldots,x_{ni}))=\bigcup\limits_{k=1}^n \downarrow(\bigvee_{i\in I}x_{ki})=\overline{\bigcup_{i\in I}\bigcup\limits_{k=1}^n\downarrow x_{ki}}=\bigvee_{i\in I}s_n(x_{1i},x_{2i},\ldots,x_{ni}).
$$
Thus $s_n$ is a Scott continuous map from $\prod\limits^nP$ into $\Gamma(P)$.

Next, let $\mathcal{U}$ be a Scott open subset of $\Gamma(P)$ and $A\in\mathcal{U}$. With loss of generality, we assume $A\not=\emptyset$. Note that since $A=\bigcup\{\da F: F\subseteq_f\!A\}$ and $\{\da F: F\subseteq_fA\}$ is  a directed family in $\Gamma(P)$, there exists a non-empty finite subset $F$ of $ A$ such that $\da \! F\in \mathcal{U}$. Let $F=\{x_1,x_2,\ldots,x_n\}$, then $s_n(x_1,x_2,\ldots,x_n)=\da\! F\in \mathcal{U}$. It follows that  $(x_1,x_2,\ldots,x_n)\in s_n^{-1}(\mathcal{U})$. Since $s_n$ is Scott continuous and $\Sigma(\prod\limits^nP)=\prod\limits^n(\Sigma P)$, there exists a family of Scott open subset $U_k, \ k=1,2,\ldots, n$ of $P$ such that
$$(x_1,x_2,\ldots,x_n)\in U_1\times U_2\times\cdots\times U_n\subseteq s_n^{-1}(\mathcal{U}).$$
Since $x_k\in A$ for $1\leq k\leq n$, we have $A\in \Diamond U_k=\{B\in \Gamma(P): \ B\cap U_k\not=\emptyset\}$. It follows that $A\in \bigcap\limits_{k=1}^n\Diamond U_k\in \upsilon(\Gamma(P))$.
For any $B\in \bigcap\limits_{k=1}^n\Diamond U_k$, there exists $y_k\in B\cap U_k$ for $1\leq k\leq n$. since $(y_1,y_2,\ldots,y_n)\in s_n^{-1}(\mathcal{U})$, we have $\bigcup\limits_{k=1}^n\da\!y_k\in \mathcal{U}$. It follows that $B\in \mathcal{U}$, i.e., $A\in\bigcap\limits_{k=1}^n\Diamond U_k\subseteq \mathcal{U}$.
$\hfill{} \square$

\begin{lem}\rm\label{corecompact}(\cite[Theorem II-4.13]{redbook})
Let $P$ be a dcpo, $(P,\sigma(P))$ is core-compact if and only if $\Sigma(Q\times P)=\Sigma Q\times\Sigma P$ for every dcpo Q.
\end{lem}

It is not difficult to see that the above lemma is still true when $P$ and $Q$ are posets.
\begin{prop}\rm\label{cor of corecompact}
Let $P$ be a poset, if $(P,\sigma(P))$ is core-compact, then $\sigma(\Gamma(P))=\upsilon(\Gamma(P))$. Moreover, $(\Gamma(P), \sigma(\Gamma(P)))$ is sober and locally compact.
\end{prop}
\noindent{\bf Proof.}
If $(P,\sigma(P))$ is core-compact, then $\Sigma(\prod\limits^nP)=\prod\limits^n(\Sigma P)$ for each $n\in\mathbb N$. It follows that $\sigma(\Gamma(P))=\upsilon(\Gamma(P))$ by Proposition \ref{Scott=upper}. Let $L=\sigma(P)$ which is a continuous lattice by asuumption and $\Gamma(P)\cong L^{op}$. By Proposition \ref{Lop is sober locally compact}, $(\Gamma(P), \sigma(\Gamma(P)))$ is a sober and locally compact space.
$\hfill{} \square$

\vskip 3mm

We have answered Problem 5.9  raised in \cite{Miao}. The following result appears in  J. Goubault-Larrecq's blog and is given by Matthew de Brecht.

\begin{lem}\rm\label{first countable product}(\cite{Matth})
Let $P$ and $Q$ be two posets. If both $(P,\sigma(P))$ and $(Q,\sigma(Q))$ are first-countable, then $\Sigma P\times\Sigma Q=\Sigma(P\times Q)$.
\end{lem}

Since a countable product of first-countable spaces is first-countable, we have the following result:

\begin{prop}\rm\label{cor of first countable}
Let $P$ be a poset, if $(P,\sigma(P))$ is first-countable, then $\sigma(\Gamma(P))=\upsilon(\Gamma(P))$.
\end{prop}
We observe that there is an adjunction between $\sigma(P)$ and $\sigma(\Gamma(P))$, which serves as a useful tool in studying the relation between $P$ and $\Gamma(P)$. For the standard theory of adjunctions or Galois connection we refer to (\cite{redbook}, Section O-3).
\begin{prop}\rm\label{retraction}
Let $P$ be a poset.
\begin{enumerate}
\item[(1)] $\eta:P\rightarrow\Gamma(P)$, $\forall x\in P, \eta(x)=\da x$. Then $\eta$ is Scott continuous.
\item[(2)] A map $f:\sigma(P)\rightarrow\sigma(\Gamma(P))$ is defined by $f(U)=\Diamond U =\{A\in\Gamma(P): A\cap U\ne\emptyset\}$. Then $f$ preserves arbitrary sups.
\item[(3)] Let the map $\eta^{-1}:\sigma(\Gamma(P))\rightarrow\sigma(P)$ be defined by $\eta^{-1}(\mathcal{U})=\{x\in P:\ \da x\in\mathcal{U}\}$. Then $\eta^{-1}\circ f=1_{\sigma(P)}$, $f\circ\eta^{-1}\leq 1_{\sigma(\Gamma(P))}$, which implys that $(\eta^{-1},f)$ is an adjunction.
\end{enumerate}
\end{prop}

\noindent{\bf Proof.}
$(1)$ Straightforword.

$(2)$ Obviously, $f$ is monotone. Let $\{U_i: i\in I\}$ be arbitrary subset of $\sigma(P)$. $f(\bigcup_{i\in I}U_i)=\{A\in\Gamma(P): A\cap\bigcup_{i\in I}U_i\neq\emptyset\}=\{A\in\Gamma(P): \exists i\in I, A\cap U_i\neq\emptyset\}=\bigvee_{i\in I}f(U_i)$.

$(3)$ Because $\eta $ is Scott continuous, $\eta^{-1}$ preserves arbitrary sups and finite infs. For any $U\in\sigma(P),
x\in\eta^{-1}(f(U))\Leftrightarrow\eta(x)\in f(U)\Leftrightarrow\da x\cap U\neq\emptyset\Leftrightarrow x\in U$, hence $\eta^{-1}\circ f=1_{\sigma(P)}$. For any $\mathcal{U}\in\sigma(\Gamma(P))$, $A\in f\circ\eta^{-1}(\mathcal U)\Leftrightarrow A\cap\eta^{-1}(\mathcal U)\neq\emptyset\Rightarrow A\in\mathcal U$, i.e., $f\circ\eta^{-1}\leq 1_{\sigma(\Gamma(P))}$. $\hfill{} \square$
\vskip 3mm
For a complete lattice $L$, the condition that $\Sigma(\prod\limits^nL)=\prod\limits^n(\Sigma L)$ for each $n\in\mathbb N$ relates to the concept of jointly continuous semilattice.
\begin{defn}\rm\label{jointly continuous}
Let $L$ be a sup semilattice, the sup operation is jointly continuous with respect to the Scott topology provided that the mapping
$$(x, y)\mapsto x\vee y: (L, \sigma(L))\times(L, \sigma(L))\rightarrow(L, \sigma(L))$$
is continuous in the product topology.
\end{defn}

\begin{prop}\rm\label{sober}(\cite[Corollary II-1.12]{redbook})
If $L$ is a dcpo and a sup semilattice
such that the sup operation is jointly Scott continuous, then $(L, \sigma(L))$ is a sober
space.
\end{prop}

Now we complete the conclusion in Proposition \ref{cor of corecompact}.
\begin{thm}\rm\label{corecom closedsubsets}
Let $P$ be a poset, then the following statements are equivalent:
\begin{enumerate}
\item[(1)] $(P,\sigma(P))$ is core-compact.
\item[(2)] $(\Gamma(P), \sigma(\Gamma(P)))$ is core-compact.
\item[(3)] $(\Gamma(P), \sigma(\Gamma(P)))$ is sober and locally compact
\item[(4)] $(\Gamma(P), \sigma(\Gamma(P)))$ is sober and locally compact with $\sigma(\Gamma(P))=\upsilon(\Gamma(P))$.
\end{enumerate}
\end{thm}

\noindent{\bf Proof.}
$(1)\Rightarrow(4)\Rightarrow(3)$. Proposition \ref{cor of corecompact}.

$(3)\Rightarrow(2)$. The open subsets of a locally compact space form a continuous lattice.

$(2)\Rightarrow(1)$. If $(\Gamma(P), \sigma(\Gamma(P)))$ is locally compact, then $\sigma(\Gamma(P))$ is a continuous lattice. By Proposition \ref{retraction}, $\sigma(P)$ is also a continuous lattice.

$(2)\Rightarrow(3)$. If $(\Gamma(P), \sigma(\Gamma(P)))$ is core-compact, then $\Sigma(\Gamma(P)\times\Gamma(P))=\Sigma\Gamma(P)\times\Sigma\Gamma(P)$ by Lemma \ref{corecompact}. Obviously, the sup operation is jointly Scott continuous. Thus $(\Gamma(P), \sigma(\Gamma(P)))$ is sober and locally compact by Proposition \ref{sober} and Lemma \ref{sober corecompact lc}.
$\hfill{} \square$
\vskip 3mm
It is proved by different method in \cite{Miao} that for a poset $P$, if $\Sigma\Gamma(P)$ is locally compact then $\Sigma L$ is core-compact. It seems that the equivalence in Theorem \ref{corecom closedsubsets} only works for posets with the Scott topology. In the following example, we give a topological space $X$ and show that  $X$ is a compact Hausdorff space while $\Sigma C(X)$ is sober but not locally compact. Moreover, $\sigma (C(X))\not=\upsilon(C(X))$.

\begin{exam}\rm\label{compact hausdorff}
Consider the subset $X=\{0\}\cup\{\frac{1}{n}: n\in\mathbb N^+\}$ of the real line $\mathbb R$, with the subspace topology. Let $C_1(X)=\{A: A\subseteq_f\{\frac{1}{n}: n\in\mathbb N^+\}\}$ and  $C_2(X)=\{B\cup\{0\}: B\subseteq\{\frac{1}{n}: n\in\mathbb N^+\}\}$, then $C(X)=C_1(X)\cup C_2(X)$. Clearly, $X$ is a compact Hausdorff space. We will prove that $C(X)$ is sober but not locally compact with respect to the Scott topology. Hence,  $\sigma (C(X))\not=\upsilon(C(X))$  also holds.

(1) We show that $\ua_{C(X)} A$ is Scott open for each $A\in C_1(X)$. Indeed, for any directed family $\{C_i\}_{i\in I}$ of closed subsets, if $\{C_i\}_{i\in I}$ is a finite family of $C_1(X)$, then there exists a largest element; if $\{C_i\}_{i\in I}$ is an infinite family of $C_1(X)$, then $\bigvee_{i\in I}C_i=\{0\}\cup\bigcup_{i\in I}C_i$; otherwise $\bigvee_{i\in I}C_i=\bigcup_{i\in I}C_i$. In all cases, if $\bigvee_{i\in I}C_i\in\ua_{C(X)} A$, there exists some $j\in I$ such that $A\subseteq C_j$.

Let $\mathcal A$ be a Scott closed subset of $C(X)$ and Max$\mathcal A$ be the set of maximal elements of $\mathcal A$. Then $\mathcal A=\da_{C(X)}\text{Max}\mathcal A$. We claim that $\mathcal A$ cannot be irreducible whenever $\text{Max}\mathcal A$ is an infinite subset of $C(X)$, which implies that $C(X)$ is sober with respect to the Scott topology.

Case 1. There exist $A_1, A_2$ in $C_1(X)\cap\text{Max}\mathcal A$ such that $A_1\setminus A_2\neq\emptyset$ and $A_2\setminus A_1\neq\emptyset$. Then $\ua_{C(X)}A_1$ and $\ua_{C(X)}A_2$ are two Scott open subsets of $C(X)$. But $\ua_{C(X)}A_1\cap\ua_{C(X)}A_2\cap\mathcal A=\emptyset$ by the maximality of $A_1$ and $A_2$.

Case 2. There is only one $A\in C_1(X)\cap\text{Max}\mathcal A$. Then must be some $B\in C_2(X)$ such that $A\setminus B\neq\emptyset$ and $B\setminus(A\cup\{0\})\neq\emptyset$. Let $x\in B\setminus(A\cup\{0\})\neq\emptyset$, then $\ua_{C(X)}A$ and $\ua_{C(X)}\{x\}$ are two Scott open subsets of $C(X)$ which intersect with $\mathcal A$. But $\ua_{C(X)}A\cap\ua_{C(X)}\{x\}\cap\mathcal A=\emptyset$ by the maximality of $A$.

Case 3. Max$\mathcal A\subseteq C_2(X)$. Let $B_1\cup\{0\}$ and $B_2\cup\{0\}$ are two different elements in Max$\mathcal A$. Then $B_1\setminus B_2\neq\emptyset$ and $B_2\setminus B_1\neq\emptyset$. For any $F_1\subseteq_f\!B_1$ and $F_2\subseteq_f\!B_2$, then $\ua_{C(X)}F_1$ and $\ua_{C(X)}F_2$ are two Scott open subsets of $C(X)$ which intersect with $\mathcal A$. Suppose that $\mathcal A$ is an irreducible closed subset. We have that $\ua_{C(X)}F_1\cap\ua_{C(X)}F_2\cap\mathcal A\neq\emptyset$, which implies that $F_1\cup F_2\cup\{0\}\in\mathcal A$. Since $F_1, F_2$ are arbitrary and $\mathcal A$ is a Scott closed subset, it follows that $B_1\cup B_2\cup\{0\}$ must be in $\mathcal A$. This contradicts with the maximality of $B_1\cup\{0\}$.

(2) Suppose that $(C(X), \sigma(C(X)))$ is locally compact. Let $\mathcal U$ be any Scott open set containing $\{0\}$ but not $\emptyset$. From the assumption, there exists a compact saturated subset $\mathcal K$ such that $\{0\}\in\text{int}_{\sigma(C(X))}\mathcal K\subseteq\mathcal K\subseteq\mathcal U$. It is easy to see that $\mathcal K=\ua_{C(X)}\text{Min}\mathcal K$ and for each $A\in\text{Min}\mathcal K$, $A=\{0\}$ or $A\subseteq_f\{\frac{1}{n}: n\in\mathbb N^+\}$. We claim that $\mathcal K$ is Scott open. Let $\{C_i\}_{i\in I}$ be any directed family of closed subsets with $\bigvee_{i\in I}C_i\in\mathcal K$, then $\bigvee_{i\in I}C_i\in\ua_{C(X)} A$ for some $A\in\text{Min}\mathcal K$. If $A=\{0\}$, then there exists some $i\in I$ such that $C_i\in\text{int}_{\sigma(C(X))}\mathcal K\subseteq\mathcal K$. Otherwise, $A\subseteq_f\{\frac{1}{n}: n\in\mathbb N^+\}$, there exists some $C_j\supseteq A$ by previous proof. Then by Lemma \ref{Scott compact open}, $\text{Min}\mathcal K$ is finite. But this is impossible, since $\mathcal K$ may contain some $A\subseteq_fX\setminus\bigcup\text{Min}\mathcal K$.

Moreover, according to the previous proof $(C(X), \upsilon(C(X)))$ is core-compact while $X$ is core-compact. By Lemma \ref{sober corecompact lc} and \ref{v is sober}, $(C(X), \upsilon(C(X)))$ is sober and locally compact. From this we get that $\sigma(C(X))\neq\upsilon(C(X))$. Indeed, we have the following result:
\end{exam}

\begin{prop}\rm
Let $Y$ be a $T_1$ space, denote by $C(Y)$ the lattice of all closed subsets of $Y$ with inclusion order, if $\sigma(C(Y))=\upsilon(C(Y))$, then $Y$ has the discrete topology.
\end{prop}

\noindent{\bf Proof.}
Suppose not, then there exists an infinite subset $A$ of $Y$ such that $\overline{A}\setminus A\neq\emptyset$. Let $x\in\overline{A}\setminus A$ and $\mathcal U=C(Y)\setminus(\{\emptyset\}\cup\{\{y\}: y\in A\})$, then $\{x\}\in\mathcal U$ and $\mathcal U\in\sigma(C(Y))$. According to the assumptions that $\sigma(C(Y))=\upsilon(C(Y))$, there exists finitely many open neighborhoods $\{U_i: 1\leq i\leq n\}$ of $x$ such that $\{x\}\in\Diamond U_1\cap\Diamond U_2\cap\cdots\cap\Diamond U_n\subseteq\mathcal U$. Let $V=U_1\cap U_2\cap\cdots\cap U_n$. Then $\{x\}\in\Diamond V\subseteq\Diamond U_1\cap\Diamond U_2\cap\cdots\cap\Diamond U_n\subseteq\mathcal U$ and $C(Y)\setminus\mathcal U\subseteq\{H\in C(Y): H\cap V=\emptyset\}$. It follows that $A\cap V=\emptyset$, and then $\overline A\cap V=\emptyset$, a contradiction.
$\hfill{} \square$
\vskip 3mm
Similar to Theorem \ref{dual property of opensubsets}, the correspondence between $P$ and $\Gamma(P)$ has been studied (see \cite{Hof,Law,ven,HoZhao}). Here we use the adjunction in Proposition \ref{retraction} to give a unified proof.
\begin{lem}\rm\label{lem}(\cite{lyv})
Let $P$ and $Q$ be  complete lattices, and let $g:P\rightarrow Q$ and $d:Q\rightarrow P$ be maps which preserve arbitrary sups and $g\circ d=1_{Q}$.
\begin{enumerate}
\item[(1)] If $P$ is completely distributive, then so is $Q$,
\item[(2)] If $P$ is hypercontinuous, then so is $Q$.
\end{enumerate}
\end{lem}

\begin{thm}\rm\label{dual property of closedsubsets}
Let $P$ be a poset,
\begin{enumerate}
\item[(1)] $P$ is continuous iff the lattice $\Gamma(P)$ is a continuous lattice.
\item[(2)] $P$ is quasicontinuous iff the lattice $\Gamma(P)$ is a quasicontinuous lattice.
\item[(3)] $P$ is algebraic iff the lattice $\Gamma(P)$ is an algebraic lattice.
\item[(4)] $P$ is quasialgebraic iff the lattice $\Gamma(P)$ is a quasialgebraic lattice.
\end{enumerate}
\end{thm}

\noindent{\bf Proof.}
(1) If $P$ is continuous, then $\sigma(P)$ is a completely distributive lattice by Theorem \ref{dual property of opensubsets}. Both $\sigma(P)$ and $\Gamma(P)$ are continuous lattices by Proposition \ref{CDlattice}. Conversely, if $\Gamma(P)$ is continuous, then $\sigma(\Gamma(P))$ is completely distributive, and so $\sigma(P)$ is a completely distributive lattice by Proposition \ref{retraction} and Lemma \ref{lem}. This shows that $P$ is a continuous poset.

(2) If $P$ is quasicontinuous, then $\sigma(P)$ is a hypercontinuous lattice by Theorem \ref{dual property of opensubsets}. It follows that $\Gamma(P)$ is a quasicontinuous lattice by Theorem \ref{property of hypercontinuous}. Conversely, if $\Gamma(P)$ is quasicontinuous, then $\sigma(\Gamma(P))$ is hypercontinuous, and so $\sigma(P)$ is a hypercontinuous lattice by Proposition \ref{retraction} and Lemma \ref{lem}. This shows that $P$ is a quasicontinous poset.

(3) For each nonempty $B\in\Gamma(P)$, it is easy to see that $B\lll B$ in $\Gamma(P)$ if and only if there exists some $k\in\text{K}(P)$ such that $B=\da k$. If $P$ is algebraic, then $A=\text{cl}(\bigcup\{\da k: k\in A\cap\text{K}(P)\})$. Thus $\Gamma(P)$ is an algebraic lattice. Conversely, if $\Gamma(P)$ is an algebraic lattice, $P$ is continuous by (1). Both $\sigma(P)$ and $\Gamma(P)$ are completely distributive. This implies that $P$ is algebraic.

(4) If $P$ is quasialgebraic, then by Theorem \ref{dual property of opensubsets}, $\sigma(P)$ is a hypercontinuous and algebraic lattice. Let $L=\sigma(P)$. Then by Theorem \ref{dual property of opensubsets} and Theorem \ref{property of hypercontinuous}, $\omega(L^{op})=\upsilon(L)=\sigma(L)$ is an algebraic lattice. This implies that $\Gamma(P)$ is a qusialgebraic lattice from Lemma \ref{priestly space}.
Conversely, if $\Gamma(P)$ is a quasialgebraic lattice, then by Lemma \ref{priestly space}, $\omega(L^{op})=\upsilon(L)$ is an algebraic lattice.
By (2) and Theorem \ref{dual property of opensubsets}, $P$ is quasicontinuous and $L$ is a hypercontinuous latiice. By Theorem \ref{property of hypercontinuous}, $\sigma(L)=\upsilon(L)$ and is a completely distibutive algebraic lattice. This together implies that $\sigma(P)$ is a hyperalgebraic lattice. Thus $P$ is a quasialgebraic poset, by Theorem \ref{dual property of opensubsets}.$\hfill{} \square$

\section{Distributive continuous lattice}
Now we return to continuous lattices to consider the condition that $\sigma(L^{op})=\upsilon(L^{op})$. At first, we list some well-known results.

\begin{fact}\rm\label{homo of sober}
If $X$ and $Y$ are both sober spaces and the open set lattice $\mathcal O(X)$ is isomorphic to $\mathcal O(Y)$, then $X$ is homeomorphic to $Y$.
\end{fact}

\begin{lem}\rm\label{finitly generated upperset}(\cite[Lemma III-5.7]{redbook})
Let $L$ be a quasicontinuous domain. If $A=\ua A$ is compact in the Scott topology, then every Scott open neighborhood $U$ of $A$ contains a finite set $F$ such that $A\subseteq\UUa F\subseteq\ua F\subseteq U$. Furthermore, $A$ is a directed intersection of all finitely generated upper sets that contain $A$ in their Scott interior.
\end{lem}

\begin{lem}\rm\label{lawson topology}(\cite{Gierz})
Let $L$ be a complete lattice.
\begin{enumerate}
\item[(1)] The Lawson topology $(L, \lambda(L))$ is compact and $T_1$.
\item[(2)] $L$ is quasicontinuous if and only if $(L, \lambda(L))$ is Hausdorff.
\item[(3)] A subset $M$ of $L$ is closed in the lower topology if and only if $M=\ua M$ and if for every ultrafilter $\mathscr F$ with $M\in\mathscr F$, liminf$\mathscr F=\bigvee\{\bigwedge F: F\in\mathscr F\}\in M$.
\item[(4)] A subset $M$ of $L$ is closed in the lower topology if and only if $M=\ua M$ and $M$ is closed in the Lawson topology.
\end{enumerate}
\end{lem}

\begin{prop}\rm\label{euqal of upsilon=sigma}
Let $L$ be a continuous lattice, the following conditions are equivalent:
\begin{enumerate}
\item[(1)] $\sigma(L^{op})=\upsilon(L^{op})$.
\item[(2)] The bi-Scott topology $\sigma_{Bi}(L)=\sigma(L)\vee\sigma(L^{op})$ is compact and Hausdorff.
\item[(3)] Let $Q(L)$ be the set of compact saturated sets of $(L, \sigma(L))$, ordered by reverse inclusion, then $\sigma(L^{op})\cong Q(L)$.
\end{enumerate}
\end{prop}

\noindent{\bf Proof.}
(1) $\Rightarrow$ (2). The bi-Scott topology $\sigma_{Bi}(L)=\sigma(L)\vee\sigma(L^{op})=\sigma(L)\vee\upsilon(L^{op})=\sigma(L)\vee\omega(L)=\lambda(L)$, and is compact and Hausdorff by Lemma \ref{lawson topology}.

(2) $\Rightarrow$ (3). The Lawson topology of $L$ is compact and Hausdorff by the above proposition. If the bi-Scott topology of $L$ is also compact and Hausdorff, then it is equal to the Lawson topology, since $\lambda(L)\subseteq\sigma_{Bi}(L)$. Let $U\in\sigma(L^{op})$ and $K_U=L\setminus U$, then $K_U$ is an upper set in $L$ and closed in the bi-Scott topology. This implies that $K_U$ is closed in the lower topology of $L$ and $K_U=\bigcap\{\ua F: F\subseteq_f\!L, K_U\subseteq\ua F\}$. We see that $\{\ua F: F\subseteq_f\!L, K_U\subseteq\ua F\}$ is a directed subset of $Q(L)$. For $K_U\subseteq\ua F_1, K_U\subseteq\ua F_2$, $K_U\subseteq\ua F_1\bigcap\ua F_2=\bigcup\{\ua(x\vee y): x\in F_1, y\in F_2\}$ which is an upper bound for $\ua F_1, \ua F_2$ in $Q(L)$. As the supremum of a directed subset of $Q(L)$, $K_U\in Q(L)$. For any $A\in Q(L)$, $A=\bigcap\{\ua F: F\subseteq_f\!L, A\subseteq\UUa F\}$ by Lemma \ref{finitly generated upperset}, then $U_A=L\setminus A\in\sigma(L^{op})$. It follows that $\sigma(L^{op})\cong Q(L)$.

(3) $\Rightarrow$ (1). From the assumption and the proof above, we can see that $Q(L)\cong\upsilon(L^{op})\cong\sigma(L^{op})$.
 By Lemma \ref{priestly space}, $\upsilon(L^{op})=\omega(L)$ is a continuous lattice, since $L$ is a continuous lattice. This implies that $\sigma(L^{op})$ is also a continuous lattice. Then $(L^{op}, \sigma(L^{op}))$ is a sober space by an argument similar to Theorem \ref{corecom closedsubsets}. By hypothesis and Lemma \ref{v is sober}, $(L^{op}, \upsilon(L^{op}))$ and $(L^{op}, \sigma(L^{op}))$ are both sober spaces with $\upsilon(L^{op})\cong\sigma(L^{op})$. It follows that $(L^{op}, \upsilon(L^{op}))$ and $(L^{op}, \sigma(L^{op}))$ are homeomorphic by Fact \ref{homo of sober}. There exists continuous functions $f: (L^{op}, \upsilon(L^{op}))\ra(L^{op}, \sigma(L^{op}))$ and $g: (L^{op}, \sigma(L^{op}))\ra(L^{op}, \upsilon(L^{op}))$ such that $f\circ g=g\circ f=1_{L^{op}}$. Obviously, $f, g$ preserve arbitrary sups and infs. We finish the proof by showing that any closed set in $(L^{op}, \sigma(L^{op}))$ is also a closed set in $(L, \omega(L))$. Let $A$ be a closed subset in $(L^{op}, \sigma(L^{op}))$ and $\mathscr F$ be an ultrafilter with $A\in\mathscr F$, then $A=\ua A$ and $f(\mathscr F)=\{f(F): F\in\mathscr F\}$ is also an ultrafilter with $f(A)\in f(\mathscr F)$. By (3) of the above lemma, we have $\text{liminf}f(\mathscr F)=\bigvee\{\bigwedge f(F): F\in\mathscr F\}\in f(A)$, and then $g(\text{liminf}f(\mathscr F))=g(\bigvee\{\bigwedge f(F): F\in\mathscr F\})=\bigvee\{\bigwedge g\circ f(F): F\in\mathscr F\}=\bigvee\{\bigwedge F: F\in\mathscr F\}=\text{liminf}\mathscr F\in g\circ f(A)=A$. Thus $A$ is closed in the lower topology of $L$. $\hfill{} \square$
\vskip 3mm
We have the following result that is similar to Proposition \ref{retraction}.
\begin{prop}\rm\label{retract of L}
Let $L$ be a distributive continuous lattice such that $\sigma(L^{op})=\upsilon(L^{op})$.
\begin{enumerate}
\item[(1)] Define $I: L\ra\sigma(L^{op})$ by $I(x)=L\setminus\ua x$ for each $x\in L$, then $I$ preserves arbitrary sups.\\
\item[(2)] For each $U\in\sigma(L^{op})$, $U=\bigcup\{L\setminus\ua F: F\subseteq_f\!L, L\setminus\ua F\subseteq U\}=L\setminus\bigcap\{\ua F: F\subseteq_f\!L, L\setminus\ua F\subseteq U\}=L\setminus K_U$, where $K_{U}=\bigcap\{\ua F: F\subseteq_f\!L, L\setminus\ua F\subseteq U\}=\bigcap\{\ua F: F\subseteq_f\!L, F\ll K_{U}\}$. Define $J: \sigma(L^{op})\rightarrow L$ by $J(U)=\bigvee^{\ua}_{F\ll K_{U}}\wedge F$. Then $J$ preserves arbitrary sups, $J\circ I=id_{L}$ and $I\circ J\leq id_{\sigma(L^{op})}$, which implies that $(J, I)$ is an adjunction and thus $J$ also preserves arbitrary infs.
\end{enumerate}
\end{prop}
\noindent{\bf Proof.}
(1) $I(\bigvee_{\alpha\in\Delta}x_{\alpha})=L\setminus\ua \bigvee_{\alpha\in\Delta}x_{\alpha}=\bigcup_{\alpha\in\Delta} L\setminus\ua x_{\alpha}=\bigcup_{\alpha\in\Delta}I(x_{\alpha})$.

(2) Obviously, for each $\mathcal U\in\sigma(L^{op})$, $K_{\mathcal U}=L\setminus\mathcal U$ is a compact saturated subset of $L$ and $J$ is a well-defined map.

Let $\{\mathcal U_{\beta}: \beta\in\nabla\}$ be an arbitrary family of open subsets in $(L^{op}, \upsilon(L^{op}))$. For each $\mathcal U_{\beta}, \beta\in\nabla$, let $L\setminus\mathcal U_{\beta}=K_{\beta}$. Let $L\setminus\bigcup_{\beta\in\nabla}\mathcal U_{\beta}=\bigcap_{\beta\in\nabla}K_{\beta}=K$. We will show that $J(\bigcup_{\beta\in\nabla}\mathcal U_{\beta})=\bigvee_{F\ll K}\wedge F=\bigvee_{\beta\in\nabla}\bigvee_{F\ll K_{\beta}}\wedge F=\bigvee_{\beta\in\nabla}J(\mathcal U_{\beta})$, i.e., $J$ preserve arbitrary sups. It is easy to see that $\bigvee_{F\ll K}\wedge F\geq\bigvee_{\beta\in\nabla}\bigvee_{F\ll K_{\beta}}\wedge F$. Indeed, if $F\ll K_{\beta}$ for some $\beta\in\nabla$, then $F\ll K$.
In the other direction, if $F\ll K=\bigcap_{\beta\in\nabla}K_{\beta}$ for some $F\subseteq_f\!L$, then there exists a finite subset $\{\beta_i: 1\leq i\leq n\}\subseteq\nabla, n\in\mathbb N$ such that $F\ll\bigcap^n_{i=1}K_{\beta_i}$. Since $\bigcap^n_{i=1}K_{\beta_i}=\bigcap^n_{i=1}\{\bigcap\ua G: G\subseteq_f\!L, G\ll K_{\beta_i}\}=\bigcap\{\bigcap^n_{i=1}\ua G_i: G_i\ll K_{\beta_i}, 1\leq i\leq n\}$, which is the intersection of a filtered family of finitely generated upper subsets, there exists $G_i\ll K_{\beta_i}, 1\leq i\leq n$ such that $\bigcap^n_{i=1}\ua G_i\subseteq\ua F$. Because $L$ is a distributive lattice, we have that $\bigcap^n_{i=1}\ua G_i=\bigcup\{\ua\bigvee^n_{i=1}x_i: x_i\in G_i, 1\leq i\leq n\}$ and $\bigvee^n_{i=1}\wedge G_i=\bigwedge\{\bigvee^n_{i=1}x_i: x_i\in G_i, 1\leq i\leq n\}$, thus $\wedge F\leq \bigvee^n_{i=1}\wedge G_i$.

For any $x\in L$, $J\circ I(x)=\bigvee_{F\ll x}\wedge F$. Since $L$ is a continuous lattice, we have $\bigvee_{F\ll x}\wedge F=x$.

For any $\mathcal U\in\sigma(L^{op})$, let $K_{\mathcal U}=L\setminus\mathcal U$. $I\circ J(\mathcal U)=L\setminus\ua\bigvee^{\ua}_{F\ll K_{\mathcal U}}\wedge F$. We claim that $I\circ J\leq 1_{\sigma(L^{op})}$. Indeed, $K_{\mathcal U}=\bigcap_{F\ll K_{\mathcal U}}\ua F\subseteq\ua\bigvee^{\ua}_{F\ll K_{\mathcal U}}\wedge F$, thus $I\circ J(\mathcal U)\subseteq L\setminus K_{\mathcal U}=\mathcal U$.
In conclusion, $(J, I)$ is an adjunction between $\sigma(L^{op})$ and $L$.$\hfill{} \square$

\begin{rem}\rm
There is a one-to-one correspondence between distributive continuous lattices (that is, continuous frames) and locally compact sober spaces in the sense of a duality of categories, namely Hofmann-Lawson duality (cf.\cite{Hof}).
Theorem \ref{corecom closedsubsets} shows that if $L$ is a distributive continuous lattice such that the hull-kernel topology of Spec$L$ is just the Scott topology, then $\sigma(L^{op})=\upsilon(L^{op})$ and $\Sigma(L^{op})$ is sober and locally compact.
There exists a complete lattice $W$ such that $\sigma(W)$ is a continuous lattice, but $W$ is not quasicontinuous (see \cite[Theorem VI-4.5.]{redbook}). By Theorem \ref{property of hypercontinuous}, $\sigma(L^{op})=\upsilon(L^{op})$ is strictly weaker that $\sigma(L)=\upsilon(L)$ for a distributive continuous lattice. We wonder whether the condition that $\sigma(L^{op})=\upsilon(L^{op})$ in Proposition \ref{Lop is sober locally compact} is necessary and whether it implies that Spec$L$ is a Scott space.
In the case of distributive algebraic lattices, we have the following result.
\end{rem}

\begin{prop}\rm\label{distributive algebraic lattice}
Let $L$ be a distributive algebraic lattice,  then the following conditions are equivalent:
\begin{enumerate}
\item[(1)] $L$ is a hyperalgebraic lattice,
\item[(2)] $\sigma(L^{op})=\upsilon(L^{op})$,
\item[(3)] The spectrum of $L$ is a quasialgebraic domain equipped with the Scott topology,
\item[(4)] $L^{op}$ is a quasialgebraic lattice.
\end{enumerate}
\end{prop}
\noindent{\bf Proof.}
(1) $\Leftrightarrow$ (3) $\Leftrightarrow$ (4). (cf. \cite{ven})

(1) implies (2). Since (1) $\Rightarrow$ (3), $L^{op}$ is isomorphic to the lattice of Scott closed subsets of a quasialgebraic domain. Then we have $\sigma(L^{op})=\upsilon(L^{op})$ by Proposition \ref{cor of corecompact}, since any quasialgebraic domain is locally compact with respect to the Scott topology.

(2) implies (1). $\sigma(L^{op})=\upsilon(L^{op})=\omega(L)$, where $\omega(L)$ is an algebraic lattice by Lemma \ref{priestly space}. So we have that $L^{op}$ is a quasialgebraic lattice from Theorem \ref{dual property of opensubsets} and that $\sigma(L^{op})$ is a hyperalgebraic lattice by Lemma \ref{approximate of hypercontinuous}. There exists an adjunction between $\sigma(L^{op})$ and both preserve arbitrary sups by Proposition \ref{retract of L}. Thus we get the conclusion that $L$ is a hyperalgebraic lattice.
$\hfill{} \square$
\vskip 3mm
Next, we construct a complete lattice $P$ such that $(P, \upsilon(P))$ is sober and locally compact, but the Scott topology on the lattice of closed subsets is not equal to the lower Vietoris topology. $P$ endowed with the Scott topology is neither locally compact nor first-countable still we have $\sigma(\Gamma(P))=\upsilon(\Gamma(P))$.

\begin{exam}\rm\label{not locally compact dcpo}
Let $P=(\mathbb N\times\mathbb N)\cup\{\top\}$ with a partial order defined as follows:
\begin{enumerate}
\item[(\romannumeral1)] $\forall(m, n)\in\mathbb N\times\mathbb N, \bot\leq(m, n)\leq\top$;
\item[(\romannumeral2)] $\forall(m_1, n_1), (m_2, n_2)\in\mathbb N\times\mathbb N, (m_1, n_1)\leq(m_2, n_2)$ iff $m_1=m_2$ and $n_1\leq n_2$.
\end{enumerate}

(1) When considering $P$ with the upper topology and let $L=\upsilon(P)$, it is easy to see that $(P, \upsilon(P))$ is a sober space. We also have that $(P, \upsilon(P))$ is locally compact, indeed $L$ is a distributive algebraic lattice. Obviously, $P$ is isolated in $L$. For any finite subset $F$ of $P$, we will show that $P\setminus\da F\ll P\setminus\da F$ in $\upsilon(P)$ by the Alexander's Lemma.
Let $\{x_i: i\in I\}$ be a subset of $P$ such that $P\setminus\da F\subseteq\bigcup_{i\in I}P\setminus\da x_i$, i.e., $\bigcap_{i\in I}\da x_i\subseteq\da F$. If $\bigcap_{i\in I}\da x_i=\{\bot\}$, then there exists $x_i, x_j, i, j\in I$ such that $\da x_i\cap\da x_j=\{\bot\}$. If $\bigcap_{i\in I}x_i\neq\{\bot\}$, then $\{x_i: i\in I\}$ has a minimal element $x_{i_0}\neq\bot$. In the first case $P\setminus\da F\subseteq P\setminus\da x_i\cup P\setminus\da x_j$ and in the second case $P\setminus\da F\subseteq P\setminus\da x_{i_0}$, which together imply that $L$ is an algebraic lattice.
Thus $\sigma(L^{op})\neq\upsilon(L^{op})$ by Proposition \ref{distributive algebraic lattice}, since $P$ is not quasicontinuous. It is not hard to verify that $\Sigma(L^{op})$ is a sober but not locally compact space.
\vskip 3mm
(2) When considering $P$ with the Scott topology, we have $\sigma(\Gamma(P))=\upsilon(\Gamma(P))$.

One can easily check that $(P, \sigma(P))$ is neither locally compact nor first-countable. For any $\mathcal U\in\sigma(\Gamma(P))$, let $\mathcal A=\Gamma(P)\setminus\mathcal U$. We will show that $\mathcal A$ is a closed subset of $(\Gamma(P), \upsilon(\Gamma(P))$. Let $A^*=\bigcup\{A: A\in\text{Max}\mathcal A\}$, then $A^*\in\Gamma(P)$ since $A^*=\eta^{-1}(\mathcal A)$. If $\top\in A^*$, then $\mathcal A=\Gamma(P)$. Otherwise,  $(\{i\}\times\mathbb N)\cap A^*$ is a finite subset for each $i\in\mathbb N$. $A^*$ with the order inherited form $P$ is an algebraic dcpo. We see that $\mathcal A$ is a closed subset of $(\Gamma(A^*), \sigma(\Gamma(A^*)))$, and is also a closed subset of $(\Gamma(A^*), \upsilon(\Gamma(A^*)))$ since $\sigma(\Gamma(A^*))=\upsilon(\Gamma(A^*))$. This means that $\mathcal A$ can be represented as the intersection of a family of finitely generated lower sets in $\Gamma(A^*)$. Thus $\mathcal A$ can be represented as the intersection of the family of finitely generated lower sets in $\Gamma(P)$, which implies that $\mathcal A$ is a closed subset of $\Gamma(P)$ with the lower Vietoris topology.
\end{exam}

Hertling \cite{Peter} constructs a complete lattice that is not jointly Scott continuous. Actually, the construction of Hertling works for any complete lattice that is not core-compact with respect to the Scott topology.

\begin{lem}\rm(\cite[Theorem II-4.10]{redbook})
Let $X$ be a $T_0$ space, the set $E=\{(x, U)\in X\times\mathcal O(X): x\in U\}$ is open in $X\times\Sigma\mathcal O(X)$ iff $X$ is core-compact.
\end{lem}

\begin{prop}\rm\label{not jointly continuous}
Let $P$ be a complete lattice and $L=P\times\sigma(P)$. If $(P, \sigma(P))$ is not core-compact, then $L$ is not jointly Scott continuous and then $\Sigma(L\times L)\neq\Sigma L\times\Sigma L$.
\end{prop}
\noindent{\bf Proof.}
By the above lemma, the set $E=\{(x, U)\in P\times\sigma(P): x\in U\}$ is not open in $\Sigma P\times\Sigma(\sigma(P))$.
Suppose that the map sup: $(L, \sigma(L))\times(L, \sigma(L))\rightarrow(L, \sigma(L)), ((x, U), (y, V))\mapsto (x\vee y, U\cup V)$ is continuous.
Since $E$ is open in $\Sigma L$, sup$^{-1}(E)$ is open in $(L, \sigma(L))\times(L, \sigma(L))$. For any $(a, W)\in E$, we have
$(a, \emptyset)\vee(\bot_P, W)=(a, W)\in E$. There exists $D_1, D_2\in\sigma(L)$ such that $((a, \emptyset), (\bot_P, W))\in D_1\times D_2\subseteq\text{sup}^{-1}(E)$.

Let $E_1=\{x\in P: (x, \emptyset)\in D_1\}$ and $E_2=\{U\in\sigma(P): (\bot_P, U)\in D_2\}$, then $a\in E_1, E_1\in\sigma(P)$ and $W\in E_2, E_2\in\sigma(\sigma(P))$. For each $b, V\in E_1\times E_2$, $(b, V)=(b, \emptyset)\vee(\bot_P, V)\in\text{sup}(D_1\times D_2)\subseteq E$. Contradiction. Thus $L$ is not jointly Scott continuous and then $\Sigma(L\times L)\neq\Sigma L\times\Sigma L$.
$\hfill{} \square$
\begin{figure}[ht]
  \centering
  \includegraphics[height=2.2in,width=2.8in]{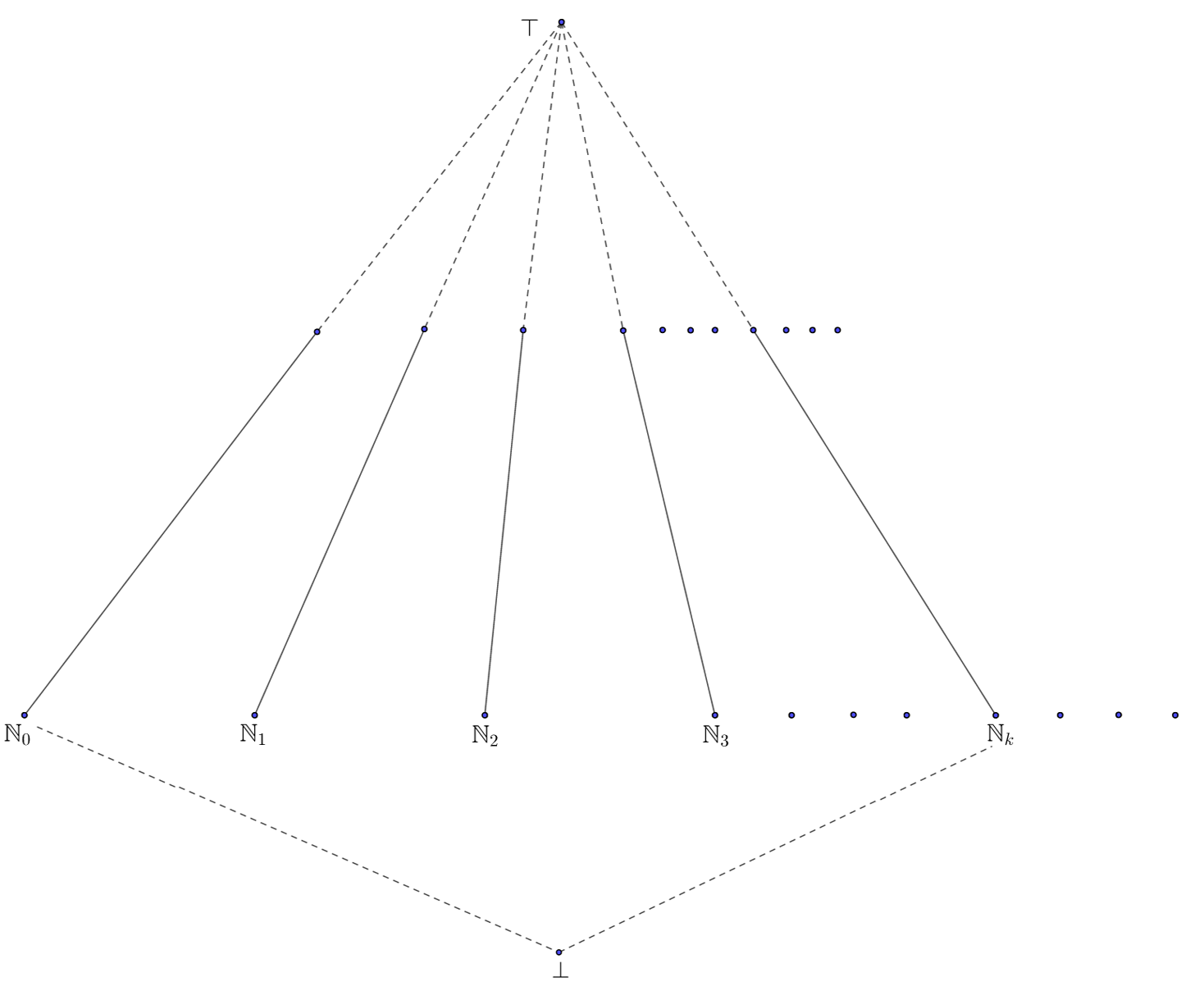} 
  \caption{Complete lattice $P$ in Example \ref{not locally compact dcpo}} 
  \label{img1} 
\end{figure}

In the following, we use the  complete lattice $P$ in Example \ref{not locally compact dcpo} to construct a  complete lattice $L$  such that $(L,\sigma(L))$ is sober but not  jointly Scott continuous.

\begin{exam}\rm
Let $P$ be the complete lattice in Example \ref{not locally compact dcpo}. Let $M=\{\bot, \top\}\cup(\mathbb N\ra\mathbb N)$ with a partial order defined as follows:
\begin{enumerate}
\item[(\romannumeral 1)] $\forall f\in(\mathbb N\ra\mathbb N), \bot\leq f\leq\top$;
\item[(\romannumeral 2)] $\forall f, g\in(\mathbb N\ra\mathbb N), f\leq g$ iff $\forall i\in\mathbb N, f(i)\geq g(i)$.
\end{enumerate}
It is easy to see that $M\cong\sigma(P)$. Let $L=P\times M$, $L$ is not jointly Scott continuous by Proposition \ref{not jointly continuous}. Let $A$ be a non-empty Scott closed subset of $L$, then $A=\da\text{Max}A$, where Max$A$ denotes the set of all maximal elements of $A$. We will show that $\Sigma L$ is a sober space by proving that $A$ can not be an irreducible closed set when $|\text{Max}A|\geq 2$.

Case 1. $|\pi_M(\text{Max}A)|=1$. $\pi_P(A)$ must be a Scott closed subset of $P$ with $|\text{Max}\pi_P(A)|\geq 2$. Since $\Sigma P$ is a sober space, there exists $B, C\in\Gamma(P)$ such that $\pi_P(A)\subseteq B\cup C$ but $\pi_P(A)\nsubseteq B$ and $\pi_P(A)\nsubseteq C$. It follows that $A\subset B\times\pi_M(A)\cup C\times\pi_M(A)$ but $A$ does not contained in any of them.

Case 2. $|\pi_M(\text{Max}A)|\geq 2$ and $\top\in\pi_M(\text{Max}A)$. Let $A_0=\{(x, y)\in\text{Max}A: y=\top\}$ and $B_0=\{(x, y)\in\text{Max}A: y\neq\top\}$, then $A_0\neq\emptyset$, $B_0\neq\emptyset$ and $A\subseteq\da A_0\cup\da\bigvee B_0$. We claim that $\da A_0$ is a Scott closed subset of $L$. Suppose not, there exists a directed subset $D\subseteq\da A_0$ such that $\bigvee D\notin\da A_0$. Let $D'=\{(\pi_P(d), \top): d\in D\}$, then $D'$ is a directed subset in $\da A_0$. We have that $\bigvee D'\notin\da A_0$ since $\bigvee D'\geq\bigvee D$, and that $\bigvee D'\notin\da\bigvee B_0$ since $\pi_M(\bigvee B_0)<\pi_M(\bigvee D')=\top$. Contradiction. Obviously, $A\nsubseteq\da A_0$ and $A\nsubseteq\da\bigvee B_0$.

Case 3. $|\pi_M(\text{Max}A)|\geq 2$ and $\top\notin\pi_M(\text{Max}A)$. Let $g=\pi_M(\bigvee A)$, then $g\in(\mathbb N\ra\mathbb N)$.

~~Case 3.1. $\bot\notin\pi_M(\text{Max}A)$. By assumption, there exists some $i_0\in\mathbb N$ such that
$B_0=\{(x, y)\in\text{Max}A: y(i_0)>g(i_0)\}\neq\emptyset$. Let $A_0=\{(x, y)\in\text{Max}A: y(i_0)=g(i_0)\}$, then $A_0\neq\emptyset$ and $A\subseteq\da A_0\cup\da\bigvee B_0$. We claim that $\da A_0$ is a Scott closed subset of $L$. Suppose not, there exists some directed subset $D\subseteq\da A_0$ such that $\bigvee D\notin\da A_0$. For each $d\in D$, define $f_d: \mathbb N\ra\mathbb N, ~f_d(n)=\pi_M(d)(n)$ if $n\neq i_0$; $f_d(n)=g(i_0)$ if $n=i_0$. Let $D'=\{(\pi_P(d), f_d): d\in D\}$, then $D'$ is a directed subset in $\da A_0$ with $\bigvee D'\geq\bigvee D$. We have that $\bigvee D'\notin\da A_0$ since $\bigvee D\notin\da A_0$,  and that $\bigvee D'\notin\da\bigvee B_0$ since $\pi_M(\bigvee B_0)<\pi_M(\bigvee D')$.
Contradiction. Obviously, $A\nsubseteq\da A_0$ and $A\nsubseteq\da\bigvee B_0$. Thus $A$ is not irreducible.

~~Case 3.2. $\bot\in\pi_M(\text{Max}A)$ and $|\pi_M(\text{Max}A)\cap(\mathbb N\ra\mathbb N)|=1$. Let $(x, g)\in\text{Max}A$, then $A\subseteq\da(x, g)\cup\pi_P(A)\times\{\bot\} $. Obviously, $A\nsubseteq\da(x, g)$ and $A\nsubseteq\pi_P(A)\times\{\bot\}$.

~Case 3.3. $\bot\in\pi_M(\text{Max}A)$ and $|\pi_M(\text{Max}A)\cap(\mathbb N\ra\mathbb N)|=1$. There exists some $i_0\in\mathbb N$ such that
$B_0=\{(x, y)\in\text{Max}A: y(i_0)>g(i_0)\}\neq\emptyset$. Let $A_0=\{(x, y)\in\text{Max}A: y(i_0)=g(i_0)\}$, then $A_0\neq\emptyset$.
Let $A_1=\da A_0\cup\pi_P(A)\times\{\bot\}$ and $B_1=\da B_0\cup\pi_P(A)\times\{\bot\}$. Similar to case 3.1, we have that $A_1$ and $B_1$ are both Scott closed subsets of $L$ with $A\subseteq A_1\cup B_1$. But $A$ is contained in non of them.

So we prove that $\Sigma L$ is a sober space while $L$ is not jointly Scott continuous. And we can see from the proof that for any dcpo $Q$, if $\Sigma Q$ is sober that $\Sigma(Q\times M)$ is also a sober space.
\end{exam}

\section{Consonance of the lower powerspace}

Given a topological space $X$, we compared the Scott topology on $C(X)$ with the lower Vietoris topology on $C(X)$. Where the latter is often called the lower powerspace over $X$ (also called the Hoare powerspace) and is denoted by $P_H(X)$. Using $Q(X)$ to denote the set of all compact saturated subsets of $X$, the upper powerspace over $X$ is $Q(X)$ with the upper Vietoris topology, which is generated by $\Box U=\{K\in Q(X): K\subseteq U\}$ as a basis, $U\in\mathcal O(X)$. $Q(X)$ with the upper Vietoris topology is also called the Smyth powerspace over $X$ and is denoted by $P_S(X)$. For each $K\in Q(X)$, let $\Phi(K)=\{U\in\mathcal O(X): K\subseteq U\}$. Then $\Phi(K)$ is a Scott open filter of $\mathcal O(X)$ and $K=\bigcap\Phi(K)$. The concept of consonance is as follows.

\begin{defn}\rm\label{consonant}
A topological space $X$ is consonant if and only if for every $\mathcal H\in\sigma(\mathcal O(X))$ and every $U\in\mathcal H$ there exists $K\in Q(X)$ such that $U\in\Phi(K)\subseteq\mathcal H$.
\end{defn}

For a sober space $X$, by the Hoffmann-Mislove theorem (see \cite[Theorem II-1.20]{redbook}),  $X$ is consonant if and only if the Scott topology on $\mathcal O(X)$ has a basis consisting of Scott open filters. Recently, Brecht and Kawai \cite{Matthew} proved that for a topological space $X$, the consonance of $X$ is equivalent to the commutativity of the upper and lower powerspaces in the sense that $P_H(P_S(X))\cong P_S(P_H(X))$ under a naturally defined homeomorphism. In that paper, they asked the following question: if $X$ is consonant, is $P_H(X)$ also consonant? We will give a partial answer below.

\begin{thm}\rm\label{consonance of PH(X)}
Let $X$ be a consonant topological space. If $\Sigma(\prod\limits^n\mathcal O(X))=\prod\limits^n(\Sigma\mathcal O(X))$ for each $n\in\mathbb N$, then $P_H(X)$ is consonant.
\end{thm}

\noindent{\bf Proof.}
Firstly, for each $n\in\mathbb N$, we define a map $\varphi_n: \prod\limits^n\mathcal O(X)\ra\upsilon(C(X))$ as follows:
$$\textstyle\forall(U_1, U_2, \ldots, U_n)\in\prod\limits^n\mathcal O(X), \displaystyle\varphi_n(U_1, U_2, \ldots, U_n)=\bigcap\limits_{k=1}^n\Diamond U_i.$$
$\varphi_n$ preserves arbitrary sups, since $\Diamond(\bigcup_{i\in I}U_i)=\bigcup_{i\in I}\Diamond U_i$ for any family $\{U_i\in\mathcal O(X): i\in I\}$. Thus $\varphi_n$ is Scott continuous and is a continuous map from $\prod\limits^n(\Sigma\mathcal O(X))$ to $\Sigma(\upsilon(C(X)))$, since $\Sigma(\prod\limits^n\mathcal O(X))=\prod\limits^n(\Sigma\mathcal O(X))$ by assumption.

Generally, $\Diamond U\cap\Diamond V\neq\Diamond(U\cap V)$ for $U, V\in\mathcal O(X)$. Let $\mathscr B=\{\bigcap_{U\in\mathcal F}\Diamond U: \mathcal F\subseteq_f\mathcal O(X)\}$, which serves as a base for $P_H(X)$. It is easy to see that $\mathscr B$ is closed for finite union.

For any $\mathscr A\in\sigma(\upsilon(C(X)))$ and $\mathcal A\in\mathscr A$. Then $\mathcal A$ is equal to the union of a family of elements in $\mathscr B$. Without loss of generality, we assume that $\mathcal A\neq\emptyset$. Here we use a trick that is frequently used.
Let $fin(\mathcal A)$ be the set of all finite unions of elements of the family, then $fin(\mathcal A)$ is a directed subset of $\mathscr B$, since $\mathscr B$ is closed for finite union. And then $\mathcal A=\bigcup fin(\mathcal A)=\bigvee^{\ua}_{\upsilon(C(X))}fin(\mathcal A)\in\mathscr A$, where $\mathscr A$ is Scott open. There exists some $\{U_j: 1\leq j\leq n\}\subseteq\mathcal O(X)\setminus\{\emptyset\}$ such that $\bigcap_{j=1}^n\Diamond U_j\in fin(\mathcal A)\bigcap\mathscr A$, i.e., $\varphi_n(U_1, U_2, \ldots, U_n)\in\mathscr A$. Since $\varphi_n$ is a continuous map from $\prod\limits^n(\Sigma\mathcal O(X))$ to $\Sigma(\upsilon(C(X)))$ by previous proof, there exists $\mathcal H_j\in\sigma(\mathcal O(X)), j=1, 2, \ldots, n$ such that
$$(U_1, U_2, \ldots, U_n)\in \mathcal H_1\times\mathcal H_2\times\cdots\times\mathcal H_n\subseteq\varphi_n^{-1}(\mathscr A).$$
By assumption, $X$ is consonant. For each $1\leq j\leq n$, there is $K_j\in Q(X)$ such that $U_j\in\Phi(K_j)\subseteq\mathcal H_j$.
Obviously, the function $\eta: x\mapsto\da x, X\ra P_H(X)$ is a topological embedding.

{\bf Claim 1}\quad $\bigcap_{j=1}^n\ua_{C(X)}\eta(K_j)$ is compact and saturated in $P_H(X)$. $\bigcap_{j=1}^n\ua_{C(X)}\eta(K_j)$ is an upper set, we only need to prove the compactness. For any $C\in C(X)$, $C\in\bigcap_{j=1}^n\ua_{C(X)}\eta(K_j)$ if and only if there exists $(x_1, x_2, \ldots, x_n)\in K_1\times K_2\times\cdots\times K_n$ such that $\da x_1\cup\da x_2\cup\cdots\cup\da x_n\subseteq C$.
The function $(x_1,\, x_2,\, \ldots,\, x_n)\mapsto\\\da x_1\cup\da x_2\cup\cdots\cup\da x_n$ from $\prod\limits^n X$ to $P_H(X)$ is continuous. Thus $\bigcap_{j=1}^n\ua_{C(X)}\eta(K_j)=\ua_{C(X)}\{\da x_1\cup\da x_2\cup\cdots\cup\da x_n: (x_1, x_2, \ldots, x_n)\in K_1\times K_2\times\cdots\times K_n\}$ is compact in $P_H(X)$ since $K_1\times K_2\times\cdots\times K_n$ is compact in $\prod\limits^n X$.

{\bf Claim 2}\quad $\bigcap_{j=1}^n\Diamond U_j\in\Phi(\bigcap_{j=1}^n\ua_{C(X)}\eta(K_j))\subseteq\mathscr A$. We have that $\{\da x_1\cup\da x_2\cup\cdots\cup\da x_n: (x_1, x_2, \ldots, x_n)\in K_1\times K_2\times\cdots\times K_n\}\subseteq\bigcap_{j=1}^n\Diamond U_j$, since $K_j\subseteq U_j$ for each $1\leq j\leq n$. Thus $\bigcap_{j=1}^n\Diamond U_j\in\Phi(\bigcap_{j=1}^n\ua_{C(X)}\eta(K_j))$, where $\Phi(\bigcap_{j=1}^n\ua_{C(X)}\eta(K_j))$ is a Scott open filter by definition. According to the previous proof, in order to prove $\Phi(\bigcap_{j=1}^n\ua_{C(X)}\eta(K_j))\subseteq\mathscr A$, we only need to prove  $\Phi(\bigcap_{j=1}^n\ua_{C(X)}\eta(K_j))\bigcap\mathscr B\subseteq\mathscr A$. Suppose that $\bigcap_{j=1}^n\ua_{C(X)}\eta(K_j)\subseteq\bigcap_{V\in\mathcal F}\Diamond V$ for some $\mathcal F\subseteq_f\mathcal O(X)$. Then for each $V\in\mathcal F$, there is some $1\leq j\leq n$ such that $K_j\subseteq V$. Otherwise, for any $1\leq j\leq n$, there is $y_j\in K_j\setminus V$. Then $\da y_1\cup\da y_2\cup\cdots\cup\da y_n\in\bigcap_{j=1}^n\ua_{C(X)}\eta(K_j)$ but $\da y_1\cup\da y_2\cup\cdots\cup\da y_n\notin\Diamond V$, which contradicts the assumption. For each $1\leq j\leq n$, let $\mathcal F_j=\{V\in\mathcal F: K_j\subseteq V\}$ and let $V_j=U_j\cap\bigcap_{V\in\mathcal F_j}V$. Then $V_j\in\Phi(K_j)\subseteq\mathcal H_j$ for each $j$, which implies that
$$(V_1, V_2, \ldots, V_n)\in \mathcal H_1\times\mathcal H_2\times\cdots\times\mathcal H_n\subseteq\varphi_n^{-1}(\mathscr A).$$
$\bigcap_{j=1}^n\Diamond V_j=\varphi_n(V_1, V_2, \ldots, V_n)\in\mathscr A$ and it is obvious that $\bigcap_{j=1}^n\Diamond V_j\subseteq\bigcap_{V\in\mathcal F}\Diamond V$. Thus $\bigcap_{V\in\mathcal F}\Diamond V\in\mathscr A$, since $\mathscr A$ is an upper set in the order of set inclusion.

Since $\bigcap_{j=1}^n\ua_{C(X)}\eta(K_j)\subseteq\bigcap_{j=1}^n\Diamond U_j\subseteq\mathcal A$, we are done.
$\hfill{} \square$
\vskip 3mm

\begin{cor}\rm
Let $X$ be a consonant topological space. If $\Sigma(\mathcal O(X))$ is core-compact or first-countable, then $P_H(X)$ is consonant.
\end{cor}

A natural question is whether consonance of $P_H(X)$ implies that $X$ is consonant. We  answer this question  negatively  by a counterexample. First, we show the following result.

\begin{prop}\rm\label{consonant corecompact}
Let $X$ be a topological space, $X$ is locally compact if and only if $X$ is core-compact and consonant.
\end{prop}

\noindent{\bf Proof.} $(\Rightarrow)$ If $X$ is a locally compact space, then $X$ is core-compact. For any $\mathcal H\in\sigma(\mathcal O(X))$ and any $U\in\mathcal H$. Suppose that $U$ is not empty, then for each $x\in U$ there is some $K_x\in Q(X)$ such that $x\in K_x^{\circ}\subseteq K_x\subseteq U$. $\{\bigcup_{x\in F}K_x^{\circ}: F\subseteq_f U\}$ is a directed subset of $\mathcal O(X)$ with $\bigvee\{\bigcup_{x\in F}K_x^{\circ}: F\subseteq_f U\}=\bigcup\{\bigcup_{x\in F}K_x^{\circ}: F\subseteq_fU\}=U\in\mathcal H$. There exists some $F\subseteq_fU$ such that $\bigcup_{x\in F}K_x^{\circ}\in\mathcal H$, since $\mathcal H$ is Scott open. Then $\bigcup_{x\in F}K_x\in Q(X)$ and $U\in\Phi(\bigcup_{x\in F}K_x)\subseteq\mathcal H$.

$(\Leftarrow)$ Let $X$ be a core-compact and consonant space. For any $x\in X$ and any open set $U$ containing $x$, there exists some open set $V$ such that $x\in V\ll U$. We employ the interpolation property to find a sequence $\{V_i\}_{i\in\mathbb N}$ of open sets such that $x\in V\ll\cdots\ll V_n\ll V_{n-1}\ll\cdots\ll V_1\ll U$ (see \cite{redbook}). Let $\mathcal H=\{W\in\mathcal O(X): W\supseteq V_i\,\text{for some} \,i\in\mathbb N\}$. It is easy to see that $\mathcal H$ is a Scott open filter with $U\in\mathcal H$. There exists some $K\in Q(X)$ such that $U\in\Phi(K)\subseteq\mathcal H$. Then $x\in V\subseteq\bigcap\mathcal H\subseteq K=\bigcap\{W\in\mathcal O(X): K\subseteq W\}\subseteq U$.
$\hfill{} \square$

\begin{exam}\rm
Exercise V-5.25 of \cite{redbook} gives  an example of a $T_0$ space $X$ such that  $\mathcal O(X)$ is a continuous lattice but $X$ itself is not  locally compact. By Proposition \ref{consonant corecompact}, $X$ is not consonant.
Let $L=\mathcal O(X)$, then $C(X)\cong L^{op}$ and $\upsilon(C(X))\cong \upsilon(L^{op})=\omega(L)$. By Lemma \ref{priestly space}, $\upsilon(C(X))$ is a continuous lattice. According to Lemma \ref{sober corecompact lc} and \ref{v is sober}, $P_H(X)$ is a sober and locally compact space. We can see from Proposition \ref{consonant corecompact} that $P_H(X)$ is consonant while $X$ is not.
\end{exam}

\vskip 3mm
\noindent{\bf Acknowledgement}
\vskip 3mm
We thank the anonymous referee for improving the presentation of the paper.

\end{document}